\documentclass{amsart} \usepackage[all]{xy}  \usepackage{amsmath,amscd, amssymb,amsfonts,pb-diagram}  \usepackage{hyperref}
\numberwithin{equation}{section}

\newcommand\st{\operatorname{st}}  \newcommand{\End}{\mbox{\rm End\,}}
\newcommand{\id}{\mathrm{id}}
\theoremstyle{plain}  
\newcommand{\G}{\mathcal{G}}  \newcommand{\C}{\mathcal{C}}
  \newcommand{\F}{\mathcal{F}}
\newcommand{\D}{\mathcal{D}} \newcommand{\M}{\mathcal{M}} 
\newcommand{\N}{\mathcal{N}}  
 \newcommand{\e}{\mathfrak{e}} 
\newtheorem{teor}{Theorem}[section] \newtheorem{corol}[teor]{Corollary}
\newtheorem{ejem-prop}[teor]{Example} \newtheorem{prop}[teor]{Proposition}
\newtheorem{lem}[teor]{Lemma}

\theoremstyle{definition} \newtheorem{defin}[teor]{Definition} 
\newtheorem{ejem}[teor]{Example}

\theoremstyle{remark} \newtheorem{obs}[teor]{Remark}

\begin{document}

\title[Clifford theory for tensor categories]{Clifford theory for tensor categories} \author{C\' esar
Galindo } \address{Departamento de  Matem\'aticas \newline \indent Pontificia Universidad
Javeriana \newline \indent Bogot\'a, Colombia.} \email{cesarneyit@gmail.com, cesar.galindo@gmail.com}   \date{\today}

\begin{abstract}
A graded tensor category over a group $G$ will be called a strongly $G$-graded tensor
category if every homogeneous component has at least one multiplicatively invertible object. Our main result is a
description of the module categories over a strongly $G$-graded tensor category as induced from module
categories over tensor subcategories associated with the subgroups of $G$.
\end{abstract}

\maketitle 
\section{Introduction}

The classical Clifford theory is an important collection of results relating representation of a group
to the representation of its normal subgroups. The principal results can be generalized using strongly
graded rings, as in \cite{dade70}.   The goal of this paper is to describe a categorical analogue of
the Clifford theory for tensor categories.

Throughout this article we work over a field $k$. By  a tensor category $(\C,\otimes,\alpha, 1 )$ we understand a $k$-linear abelian category $\C$,
endowed with a $k$-bilinear exact bifunctor $\otimes:\C\times \C\to \C$, an object $1\in \C$, and
an associativity constraint  $\alpha_{V,W,Z}: (V\otimes W)\otimes Z\to V\otimes (W\otimes Z)$, such that
Mac Lane's pentagon axiom holds \cite{BaKi},   $V\otimes 1= 1\otimes V = V$, $\alpha_{V, 1,W}= \text{id}_{V\otimes
W}$ for all $V, W\in \C$ and $\dim_k \End_\C(1)=1$.

An interesting and active problem is the classification of module categories over a tensor category.
See \cite{AM}, \cite{clas module cat on sl2}, \cite{O1}, \cite{O2}, \cite{O3}. A  left module category
over a tensor category $\C$, or a left $\C$-module category, is a $k$-linear abelian category $\M$
equipped  with an exact  bifuntor $\otimes:\C\times \M\to \M$ and natural isomorphisms
$\alpha_{X,Y,M}: (X \otimes Y)\otimes M \to X\otimes (Y\otimes M)$, $X, Y, Z \in \C, M \in \M$,
satisfying natural axioms.

\begin{defin} Let $\C$ be a tensor category, and let $\M$ be a $\C$-module category. A  $\C$-submodule
category of  $\M$ is a Serre subcategory $\N\subseteq \M$ of $\M$ such that $\N$ is a
$\C$-module category  with respect to $\otimes$.

A $\C$-module category  will be called simple if it does not contain any non-trivial $\C$-submodule
category. \end{defin} \begin{obs} A rigid tensor category over an algebraically closed field is called
finite, if it is equivalent as an abelian category to the category of finite representation of a finite
dimensional algebra, see \cite{finite-categories}. In this case the right definition of module
category is that of an \emph{exact module category}, see \textit{loc. cit.}  For exact module
categories over finite tensor categories, the notion of simple module category is equivalent to that
of indecomposable module category. In particular a semisimple module category over a fusion category
is simple if and only if it is indecomposable. \end{obs}

Let $\C$ and $\D$ be tensor categories. A $\C$-$\D$-bimodule category is a $k$-linear abelian category
$\M$, endowed with a structure of left $\C$-module category and right $\D$-module category, such that
the ``actions" commute up to natural isomorphisms in a coherent way. See Section \ref{section 2} for
details on the definitions of $\C$-module category, $\C$-bimodule category, $\C$-module functor,
$\C$-linear natural transformation and their composition.

For a right $\C$-module category $\M$ and a left $\C$-module category $\N$, the tensor
product category of $k$-linear module categories $\M\boxtimes_{\C}\N$ was defined in \cite{T-dual}; however, typically $\M\boxtimes_{\C}\N$ is not an abelian category. If
$\M$ is a $\D$-$\C$-bimodule category the category $\M\boxtimes_{\C}\N$ has a coherent left
$\D$-action.

Let $G$ be a group and $\C$ be a tensor category. We shall say that  $\C$ is  $G$-graded, if
there is a decomposition $$\C=\oplus_{x\in G}\C_x$$ of $\C$ into a direct sum of full abelian
subcategories, such that  for all $\sigma, x\in G$, the bifunctor  $\otimes$ maps $\C_\sigma\times
\C_x$ to  $\C_{\sigma x}$. See \cite{ENO}.

Recall that a graded ring $A=\oplus_{\sigma\in
G}A_\sigma$  is called strongly graded,  if  $A_xA_y= A_{xy}$ for
all $x,y\in G$. If we denote by $\C_\sigma\cdot \C_\tau\subseteq \C_{\sigma\tau}$ the
\textbf{full $k$-linear subcategory} of $\C_{\sigma\tau}$ whose objects are
direct sums of  objects of the form $V_\sigma\otimes W_\tau$, for
$V_\sigma\in \C_\sigma,$ $W_\tau\in \C_\tau$, $\sigma, \tau \in
G$, the definition of \emph{strongly graded tensor category} is
the following:
\begin{defin}\label{definicion strongly}
Let $\C= \oplus_{\sigma\in G} \C_\sigma$ be a graded tensor
category over a group $G$.  We shall say that $\C$ is \emph{strongly
graded} if the inclusion functor $\C_\sigma\cdot \C_\tau \hookrightarrow
\C_{\sigma\tau}$ is an equivalence of $k$-linear categories for all  $\sigma, \tau
\in G$.
\end{defin}
\begin{obs}
Note that $\C_\sigma\cdot \C_\tau$  is only the \emph{full $k$-linear subcategory} of $\C_{\sigma\tau}$, and not the full abelian subcategory generated by $\C_\sigma\cdot\C_\tau$. For example the  Tambara-Yamagami categories $TY(A,\chi,\epsilon)$ (see \cite{TY}) are  $\mathbb Z_2$-graded fusion categories not strongly graded. In fact,  the simple objects of $\C_0$  are invertible and $\C_1$ only have one simple. Then the objects of the category $\C_1\cdot\C_1$ has the form $(X\otimes X)^{\oplus n}$, and the full $k$-linear subcategory $\C_1\cdot\C_1$ is not equivalent to $\C_0$, if $\C_0$ has more than one simple object. Note that the abelian subcategory of $\C_0$ generated by $\C_1\cdot \C_1$ is equivalent to $\C_0$.

Also note that for every tensor category $(\C,\otimes, I)$, the $k$-linear category $\C\cdot\C$ is equivalent to $\C$, since  $V\cong I\otimes V\in \C\cdot\C$, for every  $V\in Obj(\C)$.
\end{obs}
By Lemma \ref{multi-tensorial y fuert grad}, a  graded tensor
category over a group $G$ is a strongly $G$-graded tensor
category, if and only if every homogeneous component has at least
one invertible object.
Let $\C$ be a strongly $G$-graded tensor category.  Given a $\C_e$-module category  $\M$, we shall
denote by $\Omega_{\C_e}(\M)$ the set of equivalences classes of simple $\C_e$-submodule categories of
$\M.$ By Corollary \ref{G-conjuto}, the group   $G$ acts on $\Omega_{\C_e}(\M)$ by $$G\times
\Omega_{\C_e}(\M)\to\Omega_{\C_e}(\M), \  \  (g,[X])\mapsto [\C_g\boxtimes_{\C_e} X].$$

Our main result is:

\begin{teor}[Clifford Theorem for module categories]\label{clifford} Let $\C$ be a strongly $G$-graded
tensor category and let $\M$ be a simple abelian $\C$-module category. Then:

\begin{enumerate}
  \item The action of  $G$ on $\Omega_{\C_e}(\M)$ is transitive,
  \item Let $\N$ be a simple abelian $\C_e$-submodule subcategory of $\M$. Let  $H=\st([\N])$ be
      the stabilizer subgroup of
  $[\N]\in \Omega_{\C_e}(\M)$, and let also $$\M_\N=\sum_{h\in H}\C_h\overline{\otimes} \N.$$ Then
  $\M_\N$ is a simple
  $\C_H$-module category and $\M\cong \C\boxtimes_{\C_H}\M_\N$ as $\C$-module categories.
\end{enumerate} \end{teor}

An important family of examples of strongly graded tensor categories are the crossed product tensor
categories, see \cite{nik}, \cite{tambara}. Let $\C$ be a tensor category and let $G$ be a group. We
shall denote by $\underline{G}$ the monoidal category, where the objects are the elements of  $G$,
arrows are  identities and tensor product the product of $G$.

Let $\underline{\text{Aut}_\otimes(\C)}$ be the monoidal category where objects are tensor auto-equivalences of $\C$, arrows are tensor natural isomorphisms and tensor product the composition of functors. An action of the group  $G$ over a  monoidal category  $\C$, is a
monoidal functor  $*:\underline{G}\to
\underline{\text{Aut}_\otimes(\C)}$.

Given an action $*:\underline{G}\to \underline{\text{Aut}_\otimes(\C)}$ of $G$ on $\C$,  the
$G$-crossed product tensor category, denoted by $\C\rtimes G$ is defined as follows. As an abelian
category $\C\rtimes G= \bigoplus_{\sigma\in G}\C_\sigma$, where $\C_\sigma =\C$ as an abelian
category, the tensor product is $$[X, \sigma]\otimes [Y,\tau]:= [X\otimes \sigma_*(Y), \sigma\tau],\
\  \   X,Y\in \C,\  \   \sigma,\tau\in G,$$ and the unit object is $[1,e]$. See \cite{tambara} for the
associativity constraint and a proof of the pentagon identity.

The category  $\C\rtimes G$ is $G$-graded by  $$\C\rtimes G=\bigoplus_{\sigma\in G}(\C\rtimes
G)_\sigma, \ \  \  \text{where}\  \  (\C\rtimes G)_\sigma = \C_\sigma,$$ and the objects $[1,
\sigma]\in (\C\rtimes G)_\sigma$ are invertible, with inverse $[1,\sigma^{-1}]\in (\C\rtimes
G)_{\sigma^{-1}}$.

Another useful construction of a tensor category starting from a $G$-action over a tensor category
$C$, is  the $G$-equivariantization of $\C$, denoted by $\C^G$. This construction has been used for example in
\cite{agaitsgory}, \cite{fw}, \cite{Na2}, \cite{nik}, \cite{tambara}.

The category $\C$ is a $\C\rtimes G$-module category with action $[V,\sigma]\otimes W= V\otimes \sigma_*(W)$, see \cite{nik},
\cite{tambara}. Moreover, the tensor category  of $\C\rtimes G$-linear endofunctors of $\C$ denoted by $\F_{\C\rtimes G}(\C,\C)$, is monoidally equivalent to
the $G$-equivariantization $\C^G$ of $\C$, see \cite{nik}.
With help of this equivalence can be  describe  the module categories
over $\C^G$, using the description of the module categories over the strongly $G$-graded tensor
category $\C\rtimes G$, see \cite{ENO2} for the fusion category case.

The paper is organized as follows:  Section 2 consists mainly of definitions and properties of module
and bimodule categories over tensor categories and the tensor product of module categories, that will
be need in the sequel. In Section 3 we introduce module categories graded over a $G$-set and give a
structure theorem for them.  In Section 4 the main theorem is proved. In Section 5 we describe the
simple module categories over $\C\rtimes G$ and the simple module categories over $\C^G$ if $G$ is
finite.

\section{Preliminaries}\label{section 2}

A $k$-linear category or a category additive over $k$, is a category in which the sets of arrows between two objects are $k$-vector spaces, the compositions are $k$-bilinear operations, finite direct sums exist  and there is a zero object. A $k$-linear functor $\C\to \D$ between  $k$-linear categories, is and  additive functor $k$-linear on the spaces of morphisms. The notion of $k$-bilinear bifunctor $\C\times \C'\to  \D$ is the obvious.

\begin{defin}\cite[Definition 6.]{O1} Let $\C$ be a monoidal category. A  left $\C$-module category over  $\C$, is a
category $\M$ together with a bifuntor $\otimes:\C\times \M\to \M$ and
natural isomorphisms $$ m_{X,Y,M}: (X \otimes Y)\otimes M \to X\otimes (Y\otimes M), $$ such that
\begin{align*}
    (\alpha_{X,Y,Z}\otimes M)m_{X,Y\otimes Z,M}(X\otimes m_{Y,Z,M}) &= m_{X\otimes Y,
    Z,M}m_{X,Y,Z\otimes M},\label{pentagono}\\
     1\otimes M &= M,
\end{align*} for all $X, Y, Z\in \C, M\in \M$.

A module category $\M$ over a tensor categories $\C$ always will  be abelian, and the bifunctor $\otimes:\C\times\M\to \M$ biexact. A right module category is defined in a similar way.
\end{defin}

\begin{obs}\label{obs module como funtores} For a category $\M$, the category
of $k$-linear exact endofunctors $\F(\M,\M)$ is a $k$-linear abelian  strict monoidal category, where the kerner of morphism $\tau: F\to G$ in $\F(\M,\M)$, is the functor $K: \M\to \M$, defined by $K(M)= \ker (\tau_M)$, and with the composition of functors as tensor product, and .  For a tensor category $\C$, a structure of $\C$-module  category  $(\M, \otimes, m)$ on $\M$
is the same as an exact monoidal functor $(F,\zeta): \C\to \F(\M,\M)$. The bijection is given by the
equation $V\otimes M = F(V)(M)$, identifying $$(\zeta_{V,W})_M: (F(V)\circ F(W))(M)\to F(V\otimes
W)(M)$$with $$m_{X,Y,M}^{-1}: V\otimes (W\otimes M)\to (V\otimes W)\otimes M.$$ \end{obs}
\begin{ejem}\label{example module alg} Let $(A,m,e)$ be an associative algebra in $\C$.  Let $\C_A$ be
the category of right $A$-modules  in $\C$. This is an abelian left $\C$-module category with action $
V\otimes (M,\eta)= (V\otimes M, (\text{id}_{V}\otimes \eta)\alpha_{V,M,A})$ and associativity
constraint $\alpha_{X,Y,M}$, for $X,Y\in \C, M\in \C_A$. See \cite[sec. 3.1]{O1}. \end{ejem}

\begin{ejem}
 We shall denote by  $\text{Vec}_f$ the category of finite dimensional vector spaces over $k$.
 This is a  semisimple tensor category with only one simple object. For every $k$-linear abelian
 category $\M$, there is an unique
 $\text{Vec}_f$-module category structure  with  action $k^{\oplus n}\otimes X := X^{\oplus n}$. See
 \cite[Lemma 2.2.2]{scha-tannaka}.
\end{ejem}

\begin{ejem} Let $H$ be a Hopf algebra and let $B\subseteq A$ be a left faithfully flat $H$-Galois
extension. Let $\M_B$ and  $\M^H$ be the categories of right $B$-modules and right $H$-comodules, respectively. Recall that category of right Hopf $(H,A)$-modules $\M^H_A$ is by definition the category $(\M^H)_A$ of right $A$-modules over $\M^H$. By Schneider's structure theorem \cite{schneider}, the functor $\M_B\to (\M^H)_A$,
$M\mapsto M\otimes_B A$, is a category equivalence with inverse $M\mapsto M^{\text{co}H}$. So $\M_B$
has  a $\M^H$-module category structure as in Example \ref{example module alg}. \end{ejem} For two
$\C$-modules categories  $\M$ and $\N$, a  $\C$-linear functor or module functor  $(F,\phi):\M\to \N$
consists of an exact functor  $F:\M\to \N$   and natural isomorphisms $$\phi_{X,M}:F(X\otimes M)\to
X\otimes F(M),$$ such that $$ (X\otimes \phi_{Y,M})\phi_{X,Y\otimes M}F(m_{X,Y,M})=
m_{X,Y,F(M)}\phi_{X\otimes Y, M},$$ for all $X, Y\in \C$, $M\in M$.

If $\M, \N$ are  $k$-linear abelian categories, then $\F_{\text{Vec}_f}(\M,\N)$ is the
category of $k$-linear exact functors, so $\F_{\text{Vec}_f}(\M,\N)= \F(\M,\N)$.

A $\C$-linear natural transformation  between  $\C$-linear functors $(F,\phi), (F',\phi'):\M\to \N$, is
a $k$-linear natural transformation  $\sigma:F\to F'$ such that $$\phi'_{X,M}\sigma_{X\otimes
M}=(X\otimes \sigma_M)\phi_{X,M},$$ for all $X\in \C, M\in \M$.

We shall denote the category of $\C$-linear functors and $\C$-linear natural transformations
between $\C$-modules categories $\M,  \N$ by $\F_\C(\M,\N)$. \begin{defin} Let $\C$ be a tensor
category and let $\M$ be a $\C$-module category. A  $\C$-submodule category of  $\M$ is a Serre
subcategory $\N\subseteq \M$ of $\M$, such that is a $\C$-module category  with respect to $\otimes$.

A $\C$-module category  will be called simple if it does not contain any non-trivial $\C$-submodule
category. \end{defin}

For  $\C$-linear functors $(G,\psi):\D\to \M$ and $(F,\phi):\M\to \N$,
the composition is a $\C$-linear functor $(F\circ G, \theta):\D\to \N$, where
$$\theta_{X,L}=\phi_{X,G(L)}F(\psi_{X,L}),$$ for $X\in \C$, $L\in \D$. So we have a bifunctor
\begin{gather*}
    \F_{\C}(\M,\N)\times\F_{\C}(\D,\M)\to \F_\C(\D,\N)\\
    ((F,\phi),(G,\psi))\to (F,\phi)\circ(G,\psi).
\end{gather*}

\subsection{}\label{cate modul estricta} \textbf{Strict module categories.}

A monoidal category is called \emph{strict} if its associativity constraint is the identity. In the
same way we say that a module category  $(\M,\otimes, \alpha)$ over a strict monoidal category
$(\C,\otimes ,1)$ is \emph{strict}, if $\alpha$ is the identity.

The main result of this subsection establishes that every monoidal ca\-te\-gory  $\C$ is monoidally equivalent to a strict
monoidal category $\C'$, such that every module ca\-te\-gory over $\C'$ is equivalent to a strict
one.

\begin{lem}\label{11} Let $\C$ be a monoidal category. Then  $\F_{\C}(\C,\C)\cong  \C$, where $\C$  is
a left $\C$-module category with the tensor product and the isomorphism of associativity. Moreover,
$\C^{op} \cong \F_{\C}(\C,\C)$ as monoidal categories (where $\C^{op}=\C$ as categories, and tensor product $V\otimes^{op}W = W\otimes V$). \end{lem} \begin{proof} We define the functor
$\widehat{(-)}: \C\to \F_\C(\C,\C)$ as follows: given $V\in \C$, the functor
$(\widehat{V},\alpha_{-,-V}):\C\to \C, W\mapsto W\otimes V$, $\alpha_{X,Y,V}:  \widehat{V}(X\otimes
Y)\to X\otimes \widehat{V}(Y)$ is a $\C$-module functor. If $\phi: V\to V'$ is a morphism in $\C$, we
define  the natural transformation $\widehat{\phi}: \widehat{V} \to \widehat{V'}$, as
$\widehat{\phi}_W= \text{id}_W\otimes \phi: \widehat{V}(W)=W\otimes V\to \widehat{V'}(W)=W\otimes
V'$.

The natural isomorphism $$\alpha_{-,W,V}: \widehat{V}\circ\widehat{W} \to \widehat{V\otimes^{op}W},$$
gives a structure of monoidal functor to $\widehat{(-)}$.

Let  $(F,\psi):\C\to \C$ be a module functor. Then we have a natural isomorphism
$$\sigma_X=\psi_{X,1}:F(X)=F(X\otimes 1)\to X\otimes F(1)= \widehat{F(1)}(X),$$ such that
\begin{align*}
    \alpha_{X,Y,F(1)}\sigma_{X\otimes Y} &= \alpha_{X,Y,F(1)}\psi_{X\otimes Y,1}\\
                        &= \text{id}_{X}\otimes \psi_{Y,1}\circ \psi_{X,Y}\\
                        & = \text{id}_{X}\otimes \sigma_{Y}\circ
                        \psi_{X,Y}.
\end{align*}That is, $\sigma_{X}$ is a natural isomorphism module between  $(F,\psi)$ and $(F(1),
\alpha_{-,-,F(1)})$. So the functor is essentially surjective.

Let $\phi: \widehat{V}\to \widehat{V'}$ be a $\C$-linear natural morphism. Then $\alpha_{X,1,
V}\phi_X= \text{id}_X\otimes \phi_1\alpha_{X,1,V'}$, so $\phi_X=\text{id}_X\otimes \phi_1$, and the
monoidal functor $\widehat{(-)}$ is faithful and full. Hence, by \cite[Theorem 1, p. 91]{MacLane-book}
and \cite[Proposition 4.4.2]{tannakina categories},  the functor   is an equivalence of monoidal
categories. \end{proof}

\begin{prop}\label{modulo estrict} Let $\C$ be a monoidal category, then there is a  strict monoidal
category $\overline{\C}$, such that every module category over $\overline{\C}$ is equivalent to a
strict $\overline{\C}$-module category and  $\overline{\C}$ is monoidally equivalent to $\C$.
\end{prop} \begin{proof} Let $\overline{\C}=\F_\C(\C,\C)^{op}$. By  Lemma \ref{11}, $\overline{\C}$ is
monoidally equivalent to  $\C$. Let $(\M,\otimes, m)$ be a left $\C$-module category. The category
$\F_\C(\C,\M)$ is a strict left $\overline{\C}$-module category with the composition of $\C$-module
functors. Conversely, if $\M'$ is a $\overline{\C}$-module category, then $\M'$ is a module category
over  $\C$, using the tensor equivalence $\widehat{(-}):\C\to \F_{\C}(\C,\C)$.

In a similar way to the proof of the  Lemma \ref{11}, the functor \begin{align*}
    \M\to \F_{\C}(\C,\M)\\
    M\mapsto (\widehat{M},m_{-,-,M}),
\end{align*}is an equivalence of $\C$-module categories. So every module category over $\overline{\C}$
is equiva\-lent to a strict one. \end{proof}

\subsection{Tensor product of module categories.}

\begin{defin}\cite[pp. 518]{T-dual} Let $(\M,m)$ and $(\N,n)$ be right and left $\C$-module categories
respectively. A $\C$-bilinear functor $(F,\zeta):\M\times \N\to \D$  is a  bifunctor $F:\M\times \N\to
\D,$ together with natural isomorphisms  $$\zeta_{M,X,N}:F(M\otimes X,N)\to F(M,X \otimes N),$$ such
that \begin{gather*}
    F(m_{M,X,Y},N)\zeta_{M,X\otimes Y,N}F(M,n_{X,Y,N})= \zeta_{M\otimes X,Y,N}\zeta_{M,X,Y\otimes
    N},
\end{gather*}for all $M\in \M$, $N\in \N$, $X,Y\in \C$.

A natural transformation  $\omega: (F,\zeta)\to (F',\zeta')$ between $\C$-bilinear functors, is a
natural transformation $\omega_{M,N}:F(M,N)\to F'(M,N)$ such that \begin{gather*}
   \omega_{M,X\otimes N}\zeta_{M,X,N}=\alpha_{M,X,N}'\omega_{M\otimes X,
   N},
\end{gather*}for all $M\in \M$, $N\in \N$, $X\in \C$. \end{defin}

\begin{ejem} Let $\C$ be a tensor category and let $\D$ be a tensor subcategory of $\C$. Let  $(\M,m)$
be a $\C$-module category and let $\N$ be a $\D$-module subcategory of the  $\D$-module category $\M$.
Then the functor $\C\times \N\to \M$, $(V, N)\to V\otimes M$, has a canonical  $\D$-bilinear
structure. Here, $\C$ is  a $\D$-module category in the obvious way, and the  $\D$-bilineal
isomorphism is given by  $m$. \end{ejem}

We shall denote by $\text{Bil}(\M,\N;\D)$ the category of $\C$-bilinear functors. In
\cite{T-dual}  a  $k$-linear category (not necessarily abelian) $\M\boxtimes_\C \N$ is constructed by
generators and relations, together with a $\C$-bilinear functor $T: \M\times \N\to \M\boxtimes_\C \N$,
that
 induces an equivalence of $k$-linear categories
$\F(\M\boxtimes_\C\N, \D)\to \text{Bil}(\M,\N;\D)$, for every $k$-linear category $\D$.

The objects of $\M\boxtimes_\C \N$ are finite sums of symbols  $[X ,Y]$, for  objects  $X\in
\M$, $Y\in \N$. Morphisms are sums of compositions of symbols $$[f,g]:[X,Y]\to [X',Y'],$$ for $f:X\to
X'$, $g:Y\to Y'$, symbols $$\alpha_{X,V,Y}:[X\otimes V, Y]\to [X,V\otimes Y],$$ for $X\in \M$, $V\in
\C$, $N\in \N$, and symbols for the formal inverse of $\alpha_{X,V,Y}$. The generator morphisms
satisfy the following relations:

(i) Linearity: \begin{gather*}
    [f+f',g]= [f,g]+[f',g],\  \  [f,g+g']= [f,g]+[f,g'],\\
    [a f,g]=[f,ag]=a[f,g],
\end{gather*}for all morphisms $f,f':M\to M'$ in $\M$, $g,g':N\to N'$ in $\N$, and $a\in k$.

(ii) Functoriality: \begin{gather*}
    [ff',gg']= [f',g'][f,g], \  \  [\text{id}_M,\text{id}_N]= \text{id}_{[M,N]},
\end{gather*}for all $f:M\to M'$,  $f':M'\to M''$ in $\M$, and $g:N\to N'$, $g': N'\to N''$ in $\N$.

(iii) Naturality: \begin{gather*}
    \alpha_{M',V',N'}[f\otimes u,g]=[f,u\otimes g]\alpha_{M,V,N},
\end{gather*}for morphisms $f:M\to M'$ in $\M$, $u:V\to V'$ in $\C$, and $g:N\to N'$ in $\N$.

(iv) Coherence: \begin{gather*}
    [\alpha_{M,V,W},\text{id}_N]\alpha_{M,V\otimes W, N}[\text{id}_M,\alpha_{V,W,N}]= \alpha_{M\otimes
    V,Y,N}\alpha_{M,X,Y\otimes
    N},
\end{gather*}for all $M\in \M, N\in \N, V,W\in \C$.

Let $\M,$ $\N$ be $k$-linear categories, then  the category   $\M\boxtimes \N:=
\M\boxtimes_{\text{Vec}_f}\N$, is the tensor product of  $k$-linear tensor categories; see
\cite[Definition 1.1.15]{BaKi}. If  $\M$ and $\N$ are semisimple categories, this is the Deligne's
tensor  product of abelian categories \cite{D}. \begin{defin}\cite[pp. 517]{T-dual} Let $\C_1$ and
$\C_2$ be tensor categories. A $\C_1$-$\C_2$-bimodule category is a $k$-linear abelian category $\M$,
equipped with exact bifunctors $\otimes :\C_1\times \M\to \M$, $\otimes :\M\times \C_2\to \M$,  and
naturals isomorphisms \begin{gather*}
    \alpha_{X,Y,M}: (X \otimes Y)\otimes M \to X\otimes (Y\otimes M),\\
     \alpha_{X,M,S}: (X \otimes M)\otimes S \to X\otimes (M\otimes S),\\
      \alpha_{M,S,T}: (M \otimes S)\otimes T \to M\otimes (S\otimes
      T),
\end{gather*}for all $X, Y\in \C_1$, $M\in \M$, $S, T\in \C_2$, such that $\M$ is a left $\C_1$-module
category with $\alpha_{X,Y,M}$, it is a right  $\C_2$-module category with $\alpha_{M,S,T}$, and
\begin{gather*}
    \text{id}_X\otimes \alpha_{Y,M,S}\alpha_{X,Y\otimes M, S}\alpha_{X,Y,M}\otimes \text{id}_S
 = \alpha_{X, Y, M\otimes Z}\alpha_{X\otimes Y,M,S},\\
 \text{id}_X\otimes \alpha_{M,S,T}\alpha_{X,M\otimes S, T}\alpha_{X,M,S}\otimes \text{id}_T
 = \alpha_{X, M, S\otimes T}\alpha_{X\otimes M,S,T}.
\end{gather*} \end{defin}

%\begin{obs}\label{strict bimodule} %If $\C_1$ and $\C_2$ are tensor categories then the category $\C_1\boxtimes\C_2$ is a tensor category with tensor product $$\otimes \times \otimes :(\C_1\boxtimes\C_2)\times (\C_1\boxtimes\C_2)\to \C_1\boxtimes\C_2$$and unit object $1\boxtimes 1$, is not difficult to see that a $\C_1$-$\C_2$-bimodule category is the same as a right $\C_1\boxtimes\C_2^{op}$-module category. % \end{obs}

If $\M$ is a  $(\C_1,\C_2)$-bimodule category and $\N$ is a right $\C_2$-bimodule category, then the
category $\M\boxtimes_{\C_2}\N$ has a structure of left  $\C_1$-module category.

The action of an object $X\in \C_1$ over an object  $[M,N]\in \M\boxtimes_{\C_2}\N$ is given by
$$X\otimes [M,N]= [X\otimes M, N].$$The action over the morphisms $\alpha_{M,Y,N}$ is given by
$\text{id}_X\otimes \alpha_{M,Y,N} = \alpha_{X\otimes M,X,Y,N}\circ [\alpha_{X,M,Y}^{-1},N]$, and the
associativity is $$[\alpha_{X,Y,M},N]:[(X\otimes Y)\otimes M, N]\to [X\otimes (Y\otimes M),N].$$

\begin{prop}\label{prodiedades tensor propuct module categories} Let $\C$ be a tensor category. Let
$\M_1, \M_2$ be  $\C$-bimodule categories, and let $\M_3$ be a right $\C$-module category. Then
\begin{enumerate}
  \item $\C\boxtimes_{\C}\M_3\cong \M_3$, as left $\C$-module categories.
  \item $(\M_1\boxtimes_{\C}\M_2)\boxtimes_{\C}\M_3 \cong
      \M_1\boxtimes_{\C}(\M_2\boxtimes_{\C}\M_3)$, as left $\C$-module categories.
  \item if $\M =\oplus_{i}^n\M^i$, $\N=\oplus_{j}^m\N^j$, as  right and left $\C$-module
      categories,
  then $\M\boxtimes_{\C}\N =\oplus_{i,j}\M_{i}\boxtimes_{\C}\N_j$, as  $k$-linear categories.
\end{enumerate} \end{prop} \begin{proof} By Proposition \ref{modulo estrict}, we can suppose that
all module categories are strict.

(1) The functor  $F :\M\to \C\boxtimes_{\C}\M$ $M \mapsto [1,M]$, is a category equivalence. In
effect, using the isomorphism $\alpha_{1,X,M}$, we can see that $F$ is essentially surjective, and
every morphism between $[1,M]$ and $[1,N]$ is of the form $[1,f]$, for $f:M\to N$. Then $F$ is
faithful and full. Moreover, with the natural isomorphism $\eta_{V,M}=\alpha_{1,V,M}: F(V\otimes N)\to
V\otimes F(N)$, the pair $(F,\eta)$ is a $\C$-linear functor, since \begin{align*}
    \eta_{V\otimes W,M} = \alpha_{1,V\otimes W, M} =& \alpha_{V,W,M}\circ \alpha_{1,V,W\otimes M}\\
                                                  =& \text{id}_V\otimes \alpha_{1,W,M}\circ
                                                  \eta_{V,W\otimes M}\\
                                                  =& \text{id}_V\otimes \eta_{W,M}\circ
                                                  \eta_{V,W\otimes M}.\end{align*}

(2)\  For every object $M_1\in \M_1$, the functor $\lambda_{M_1}: \M_2\times \M_3 \to
(\M_1\boxtimes_{\C}\M_2)\boxtimes_{\C}\M_3$, where \begin{equation*}
    \lambda_{M_1} (M_2, M_3) = [[M_1,M_2],M_3], \quad
    \lambda_{M_1} (f, g) = [[\text{id}_{M_1},f],g],
\end{equation*} with the natural transformation $\eta^1_{M_2, V,M_3}: =\alpha_{[M_1,M_2],V,M_3}$, is a
$\C$-bilinear functor. So we have a family of functors $\overline{\lambda_{M_1}}: \M_2\boxtimes_{\C}
\M_3 \to (\M_1\boxtimes_{\C}\M_2)\boxtimes_{\C}\M_3$, $\overline{\lambda_{M_1}} ([M_2, M_3]) =
[[M_1,M_2],M_3]$. Now, the functor \begin{align*}
    \M_1\times(\M_2\boxtimes_{\C}\M_3)&\to (\M_1\boxtimes_{\C}\M_2)\boxtimes_{\C}\M_3,\\
    (M_1,[M_2,M_3]) &\mapsto \overline{\lambda_{M_1}}([M_2,M_3]),
\end{align*} with the natural transformation $\eta^2_{M_1,V, [M_2,M_3]}= \alpha_{M_1, V, [M_2,M_3]}$,
is a $\C$-bilinear functor. So we have a functor $\pi: \M_1\boxtimes_{\C}(\M_2\boxtimes_{\C}\M_3)\to
(\M_1\boxtimes_{\C}\M_2)\boxtimes_{\C}\M_3$, $[M_1,[M_2,M_3]] \mapsto [[M_1,M_2],M_3]$. The functor
$\pi$ is essentially surjective and \begin{align*}
    \pi([f,[g,h]]) &= [[f,g],h]\\
    \pi([\text{id}_{M_1}\alpha_{M_2,V,M_3}]) &= \alpha_{[M_1,M_2],V,M_3}\\
    \pi( \alpha_{M_1,V,[M_2,M_3]}) &=
    [\alpha_{M_1,V,M_2},\text{id}_{M_3}].
\end{align*} So  $\pi$ is  faithful and full, hence  by \cite[Theorem 1, pp. 91]{MacLane-book}, the
functor  $\pi$ is a category equivalence. Finally, note that the functor $\pi$ is $\C$-linear.

(3) Its follows directly by the construction of  $\M\boxtimes_{\C}\N$. \end{proof} Let $\C$
be a $G$-graded tensor category. Note that if $H\subseteq G$ is a subgroup of  $G$, then the category
$\C_H=\oplus_{\tau \in H}\C_\tau$  is a tensor subcategory of  $\C$.

We shall say that an object $U\in \C$ is \emph{invertible} if the functor $U\otimes (-):\C\to \C,
V\mapsto U\otimes V$ is a category equivalence or, equivalently, if there is an object $U^* \in \C$,
such that $U^*\otimes U\cong U\otimes U^*\cong 1$.

\begin{prop}Let $\C$ be a $G$-graded category and let
$H\subseteq G$ be a subgroup of $G$. Suppose that the every category
$\C_{\sigma}$ has at least one
invertible object, for every $\sigma \in G$. Let
$\M$ be a module category over $\C_H= \oplus_{h\in H}\C_h$. Then
the $k$-linear category $\C\boxtimes_{\C_H}\M$ is an abelian
category. Moreover, since $\C$ is a $\C$-$\C_H$-bimodule category then $\C\boxtimes_{\C_H}\M$ is a left module category over the tensor category $\C$.\end{prop}

\begin{proof} We shall suppose that the tensor category $ \C$ is strict. Let $\Sigma= \{e,\sigma_1,\ldots\}$ be a set of representatives of the cosets $G/H$.
Since  $\C =\bigoplus_{\sigma\in \Sigma}\C_{\sigma H}$ as right $\C_H$-module categories,
$\C\boxtimes_{\C_H}\M = \bigoplus_{\sigma\in \Sigma}\C_{\sigma H}\boxtimes_{\C_H}\M$, as $k$-linear
categories, by Proposition \ref{prodiedades tensor propuct module categories}.

For every coset $\sigma H$ in $G$, let $U_\sigma\in \C_{\sigma}$ be an invertible object. The
functor $U_\sigma:\C_H\to \C_{\sigma H}$, $V\mapsto U_\sigma\otimes V$ is a category equivalence with
a quasi-inverse $U_{\sigma}^*:\C_{\sigma H}\to \C_H$, $W\to U_{\sigma}^*\otimes W$. Then we can
assume, up to isomorphisms, that every object of $\C_{\sigma H}$ is of the form $U_\sigma\otimes V$,
where $V\in \C_H$.

Let $\bigoplus_{i}[V_i, M_i]\in \C_{\sigma H}\boxtimes_{\C_{H}}\M$. For every  $V_i$ there exist $V_i'$ such that $V_i\cong U_\sigma\otimes V_i'$. Then $\bigoplus_{i}[V_i, M_i]\cong [U_\sigma, \bigoplus_{i} V_i'\otimes M_i]$, \textit{i.e.}, we can assume, up to isomorphisms, that every object
of $ \C_{\sigma H}\boxtimes_{\C_{H}} \M$ is of the form $[U_\sigma, M]$.

If $U_\sigma\otimes V\cong U_\sigma$ then $V\cong 1$; so every
morphism $[U_\sigma,M]\to [U_\sigma,M']$  is of the form
$[\text{id}_{U_\sigma},f]$, where $f:M\to M'$. Then the functor $:\M\to \C_{\sigma
H}\boxtimes_{\C_{H}}\M, f\mapsto [\text{id}_{U_\sigma},f]$ is an equivalence of $k$-linear categories. We define the abelian structure over $\C_{\sigma
H}\boxtimes_{\C_{H}}\M$ as the induced by this equivalence.

For the second part, note that $$\C\boxtimes_{\C_H}\M = \bigoplus_{\sigma\in
\Sigma}\C_{\sigma H}\boxtimes_{\C_H}\M$$ as abelian category, so we need to prove that if
\begin{equation}\label{1 suc exact}
0\to [U_\sigma,S]\to [U_\sigma,T]\to [U_\sigma,W]\to 0
\end{equation} is an exact sequence in $\C_{\sigma H}\boxtimes_{\C_{H}} \M$, then  the sequence
\begin{equation}\label{2 suc exact}
0\to [X\otimes U_\sigma,S]\to [X\otimes U_\sigma,T]\to [X\otimes U_\sigma,W]\to 0\end{equation} is exact for all $X\in \C$. Since $\C= \bigoplus_{\sigma \in G}\C_\sigma$ we can suppose that $X\in \C_\tau$, then $[X\otimes U_\sigma,S], [X\otimes U_\sigma,T], [X\otimes U_\sigma,W]\in \C_{\tau\sigma H}\boxtimes_{\C_{H}} \M$.

Let $U_{\tau\sigma} \in \C_{\tau\sigma}$ with inverse object $U^*_{\tau\sigma}\in \C_{(\tau\sigma)^{-1}}$, so we have  the following commutative diagram
%$$
%\begin{diagram}
%\node{0\to[X\otimes U_\sigma,S]}\arrow{s}\arrow{e,t}{[\id,f]}\node{[X\otimes U_\sigma,T]} \arrow{s}\arrow{e,t}{[\id,g]} \node{[X\otimes U_\sigma,W]\to 0}\arrow{s}\\
%\node{0\to [U_{\tau\sigma}\otimes U_{\tau\sigma}^*\otimes X\otimes U_\sigma,S]} \arrow{s,l}{\alpha_{U_{\tau\sigma},U_{\tau\sigma}^*\otimes X\otimes U_\sigma,S}}\arrow{e,t}{[\id ,f]}\node{[U_{\tau\sigma}\otimes U_{\tau\sigma}^*\otimes X\otimes U_\sigma,T]} \arrow{s,l}{\alpha_{U_{\tau\sigma},U_{\tau\sigma}^*\otimes X\otimes U_\sigma,T}}\arrow{e,t}{[\id,g]} \node{[U_{\tau\sigma}\otimes U_{\tau\sigma}^*\otimes X\otimes U_\sigma,W]\to 0}\arrow{s,l}{\alpha_{U_{\tau\sigma},U_{\tau\sigma}^*\otimes X\otimes U_\sigma,W}}\\
%\node{0\to [U_{\tau\sigma}, (U_{\tau\sigma}^*\otimes X\otimes U_\sigma)\otimes S]}\arrow{e,t}{[\id,\id_{U_{\tau\sigma}^*\otimes X\otimes U_\sigma}\otimes f]}\node{[U_{\tau\sigma}, (U_{\tau\sigma}^*\otimes X\otimes U_\sigma)\otimes T]} \arrow{e,t}{[\id,\id_{U_{\tau\sigma}^*\otimes X\otimes U_\sigma}\otimes g]} \node{[U_{\tau\sigma}, (U_{\tau\sigma}^*\otimes X\otimes U_\sigma)\otimes W]\to 0}
%\end{diagram}
%$$

$$
\begin{diagram}
\node{0\to[X U_\sigma,S]}\arrow{s}\arrow{e,t}{[\id,f]}\node{[X U_\sigma,T]} \arrow{s}\arrow{e,t}{[\id,g]} \node{[X U_\sigma,W]\to 0}\arrow{s}\\
\node{0\to [U_{\tau\sigma} (U_{\tau\sigma}^* X U_\sigma),S]} \arrow{s,l}{\alpha_{U_{\tau\sigma},U_{\tau\sigma}^* X U_\sigma,S}}\arrow{e,t}{[\id ,f]}\node{[U_{\tau\sigma} (U_{\tau\sigma}^* X U_\sigma),T]} \arrow{s,l}{\alpha_{U_{\tau\sigma},U_{\tau\sigma}^* X U_\sigma,T}}\arrow{e,t}{[\id,g]} \node{[U_{\tau\sigma} (U_{\tau\sigma}^* X U_\sigma),W]\to 0}\arrow{s,l}{\alpha_{U_{\tau\sigma},U_{\tau\sigma}^* X U_\sigma,W}}\\
\node{0\to [U_{\tau\sigma}, (U_{\tau\sigma}^* X U_\sigma) S]}\arrow{e,t}{[\id,\id_{U_{\tau\sigma}^* X U_\sigma} f]}\node{[U_{\tau\sigma}, (U_{\tau\sigma}^* X U_\sigma) T]} \arrow{e,t}{[\id,\id_{U_{\tau\sigma}^* X U_\sigma} g]} \node{[U_{\tau\sigma}, (U_{\tau\sigma}^* X U_\sigma) W]\to 0}
\end{diagram}
$$ where tensor symbols between objects and morphism have been omitted as a space-saving measure.

Then the sequence \eqref{2 suc exact} is exact if and only if the sequence \begin{equation}\label{3 suc exact}  0\to [U_{\tau\sigma}, (U_{\tau\sigma}^* X U_\sigma) S] \to [U_{\tau\sigma}, (U_{\tau\sigma}^* X U_\sigma) T]\to  [U_{\tau\sigma}, (U_{\tau\sigma}^* X U_\sigma) W]\to 0\end{equation}is exact. By definition the sequence \eqref{3 suc exact} is exact if and only if the sequence  $$0\to  (U_{\tau\sigma}^* X U_\sigma) S \to (U_{\tau\sigma}^* X U_\sigma) T\to  (U_{\tau\sigma}^* X U_\sigma) W\to 0$$ in $\M$ is exact, but since $\M$ is a $\C_H$-module category it is exact.
\end{proof}

\section{Strongly graded tensor categories}

Recall from Definition \ref{definicion strongly} that the $G$-graded category $\C$ is called strongly
graded if the inclusion functor $\C_\sigma\cdot \C_\tau \to \C_{\sigma\tau}$ is a category
equivalence for all $\sigma, \tau \in G$.

\begin{lem}\label{multi-tensorial y fuert grad} Let  $\C$ be a tensor category. Then $\C$ is strongly
graded over $G$ if and only if  the category  $\C_\sigma$  has at least one multiplicatively invertible element,   for all
$\sigma\in G$. Moreover, in this case the Grothendieck ring of $\C$ is a $G$-crossed product.
\end{lem} \begin{proof}  If $\C$ is strongly graded by definition there  there exist objects
$V_1,\ldots, V_n \in \C_\sigma, W_1,\ldots, W_t \in \C_{\sigma^{-1}}$,   such that $1\cong
\bigoplus_{i,j}V_i\otimes W_j$, then $\End_{\C}( \bigoplus_{i,j}V_i\otimes W_j)\cong \End_{\C}(1)\cong
k$, so $n=1, t=1$. That is, there exist objects $V\in \C_\sigma,  W\in \C_{\sigma^{-1}}$, such that
$V\otimes W\cong 1$.

Conversely, suppose that $\C_\sigma$ has at least an invertible object for all $\sigma\in G$.   Let
$U_\sigma\in \C_\sigma$ be an invertible object with dual object $U_\sigma^*\in \C_{\sigma^{-1}}$, so
$V\cong U_\sigma\otimes (U_\sigma^*\otimes V)$ for every $V\in \C_{\sigma\tau}$. Then the inclusion
functor is essentially surjective, and therefore it is an equivalence.

Recall that by definition a graded ring $A=\oplus_{\sigma\in G}A_\sigma$ is a crossed product over $G$
if  for all $\sigma \in G$ the abelian group $A_\sigma$ has at least an invertible element. Thus, by
the first part of the lemma, the Grothendieck ring of $\C$ is a $G$-crossed product if $\C$ is
strongly graded. \end{proof}

\begin{ejem}\label{pointed categories} Let $\text{Vec}^G_{\omega}$ be the semisimple category of finite
dimensional $G$-graded vector spaces,  with constraint of associativity $\omega(a,b,c)\text{id}_{abc}$ for
all $a,b,c\in G$, where $\omega\in Z^3(G,k^*)$ is a 3-cocycle. Then $\text{Vec}^G_{\omega}$ is a strongly
$G$-graded tensor category. \end{ejem}

\begin{ejem}
Let $\C\rtimes G$ a crossed product tensor category. As we saw in the Introduction, the category $C\rtimes G$ is a strong $G$-graded tensor category. Now if we take a normalized 3-cocycle $ \beta \in
Z^3(G,k^*) $  and we define a new associator $\alpha^\beta_{[U,\sigma],[V,\tau],[W,\rho]}= \beta(\sigma,\tau,\rho)\alpha_{[U,\sigma],[V,\tau],[W,\rho]}$, then the new tensor category is strongly $G$-graded too.

\end{ejem}
\subsection{} \textbf{\textbf{Module categories graded over a $G$-set.}}

\begin{defin} Let $\C=\oplus_{\sigma\in G}\C_\sigma$ be a graded tensor category and let $X$ be a left
$G$-set. A  \emph{left $X$-graded $\C$-module category}  is a left $\C$-module category $\M$ endowed
with a decomposition  $$\M=\oplus_{x\in X}\M_x,$$ into a direct sum of full abelian subcategories,
such that for all $\sigma\in G$, $x\in X$, the bifunctor  $\otimes$ maps $\C_\sigma\times \M_x$ to
$\M_{\sigma x}$. \end{defin} An  \emph{$X$-graded $\C$-module functor}  $F:\M\to \N$ is a $\C$-module
functor such that $F(\M_x)$ is mapped to $\N_x$, for all $x\in X$.

\begin{defin} A left $X$-graded $\C$-submodule category of  $\M$  is serre subcategory $\N$
of $\M$ such that $\N$ is an $X$-graded $\C$-module category  with respect to $\otimes$, and the
grading $\N_x\subseteq  \M_x$, $x\in X$.

An $X$-graded $\C$-module category  will be called \emph{simple} if it  contains no nontrivial
$X$-graded $\C$-submodule category. \end{defin}

\begin{lem} Let $\C$ be a $G$-graded tensor category and  let $H\subseteq G$ a subgroup of $G$. If
$\N$ is a left $\C_H$-module category, then the category $\C\boxtimes_{\C_H}\N$ is a $G/H$-graded
$\C$-module category with grading $( \C\boxtimes_{\C_H}\N)_{\sigma H}= (\oplus_{\tau\in \sigma
H}\C_\tau) \boxtimes_{\C_H}\N$. \end{lem}

\begin{proof} Let $\Sigma= \{e,\sigma_1,\ldots\}$ be a set
of representatives of the cosets of $G$ modulo $H$. By Proposition \ref{prodiedades tensor propuct
module categories},   $\C\boxtimes_{\C_H}\N = \bigoplus_{\sigma\in \Sigma}\C_{\sigma
H}\boxtimes_{\C_H}\N$ as $k$-linear categories, and by the definition of the action of $\C$,  the
module category $\C\boxtimes_{\C_H}\N$ is $G/H$-graded. \end{proof} \begin{prop}\label{prod tensor
algebras} Let $\C$ be a strongly $G$-graded tensor category, and let $(A,m,e)$ be an algebra in
$\C_H$.  Then $\C\boxtimes_{\C_H}(\C_H)_A\cong \C_A$ as  $G/H$-graded $\C$-module categories.
\end{prop} \begin{proof}Let $\Sigma= \{e,\sigma_1,\ldots \}$ a set of representatives of the cosets of
$G$ modulo $H$. The $\C$-module category  $\C_A$ has a canonical $G/H$-grading: if $(M,\rho)$ is an
$A$-module then $$(M,\rho)=\bigoplus_{\sigma\ \in \Sigma}(M_{\sigma H},\rho_{\sigma H}),$$where
$M_{\sigma H}=\bigoplus_{h\in H}\M_{\sigma h}$, $\rho_{\sigma H}=\bigoplus_{h\in H}\rho_{\sigma h}.$

Let us consider the canonical  $\C$-linear functor $F: \C\boxtimes_{\C_H}(\C_H)_A \to \C_A$,
\begin{equation*}
 [V,(M, \rho)]\mapsto (V\otimes M,\text{id}_{V}\otimes \rho).
\end{equation*} We shall first show that $F$ is a category equivalence.

Let $U_\sigma\in \C_{\sigma H}$ be an invertible object for every coset of $H$ on $G$. Let
$(M,\rho)\in \C_A$ be a homogeneous $A$-module of degree $\sigma^{-1}H$. Then the $A$-module
$(U_\sigma\otimes M,\text{id}_{U_\sigma}\otimes\rho)$ is also an $A$-module in $\C_H$ and
$F([U_{\sigma^{-1}},( U_\sigma\otimes M,\text{id}_{U_\sigma}\otimes\rho)])\cong (M,\rho)\in \C_A$. So
$F$ is an essentially  surjective functor.

 We can suppose, up to isomorphisms, that every object of $ \C_{\sigma H}\boxtimes_{\C_{H}}
(\C_H)_A$ is of the form $[U_\sigma,(M,\rho)]$. Then $F([U_g,(M,\rho)])= (U_g\otimes
M,\text{id}_{U_g}\otimes\rho)$. Now it is clear that the functor $F$ is faithful and full, so by
\cite[Theorem 1, p. 91]{MacLane-book} the functor  $F$ is a category equivalence. \end{proof}

\begin{teor}\label{estruc simples grad} Let $\C$ be a strongly graded tensor category over a group $G$
and let $X$ be a transitive $G$-set. Let $\M$ and $\N$ be  non zero $X$-graded modules categories.
Then \begin{enumerate}
  \item $\M\cong \C\boxtimes_{\C_{H}}\M_x$ as $X$-graded $\C$-module categories, where, for all $x\in
      X$, $H=\st(x)$ is the stabilizer
subgroup of $x\in X$.
  \item There is a bijective correspondence between isomorphisms classes of $X$-graded $\C$-module
      functors $(F,\eta):\M\to \N$ and $\C_H$-module functors $(T,\rho): \M_x\to \N_x$.
\end{enumerate} \end{teor}

\begin{proof} (1) Choose  $x\in X$, and denote  $H=\st(x)$. In a similar way to the proof of
Proposition \ref{prod tensor algebras}, the canonical functor $\mu: \C\boxtimes_{\C_{H}}\M_x \to \M$
\begin{equation*}[V,M] \to V\otimes M, \end{equation*} is a category equivalence and it  respects the
grading.

The proof of part (1) of the theorem is completed by showing that the functor $\mu$ is a $\C$-module
functor. Indeed, by Proposition \ref{modulo estrict} we can assume that the module categories are
strict, hence
 \begin{align*}
    \mu(V\otimes [W,M_x]) &= \mu([V\otimes W,M_x])\\
                                         &= (V \otimes W)\otimes M_x = V\otimes(W\otimes M_x)\\
                                         &=V\otimes\mu([W, M_x]),
\end{align*} \textit{i.e.}, $\mu$ is a  $\C$-module functor.

 (2) By the first part we can suppose  $\N= \C\boxtimes _H \N_x$. Let $(F,\mu ): \N_x\to
\M_x$ be a $\C_H$-module functor, the functor \begin{align*}
    I(F):\C\times\N_x&\to \M\\
    (S, N) &\mapsto S\otimes F(N)
\end{align*} with the natural transformation $\text{id}_S\otimes\mu_{V,N}: I(F)(S, V\otimes N)\to
I(S\otimes V, N)$ is a $\C_H$-bilinear functor,  so we have a functor \begin{align*}
    I(F):\C\boxtimes_{\C_H}\N_x&\to M\\
    [S, N]&\mapsto V\otimes F(N)\\
    \alpha_{S,V,N}&\mapsto \text{id}_S\otimes \mu_{V,S},
\end{align*}and this is an $X$-graded $\C$-module functor in the obvious way.

Let $(F=\oplus_{s\in X}F_s,\eta):\C\boxtimes_{\C_H}\N_x\to M$ be an $X$-graded $\C$-module functor.
Consider the natural isomorphism $$\sigma_{[V,N]}:= \eta_{V,[1,N]}: F([V,N])\to V\otimes F_x([1,N]) =
I(F_x)([V,N]),$$ \begin{align*}
    \sigma_{X\otimes[V, N]}= \eta_{X\otimes V,[1,N]} &= \text{id}_X\otimes \eta_{V,[1,N]}\circ
    \eta_{X,[V,N]}\\
                                                     &= \text{id}_X\otimes \sigma_{[V,N]}\circ
                                                     \eta_{X,[V,N]}.
\end{align*} So $\sigma$ is a natural isomorphism of module functors.\end{proof}

\begin{corol} Let $\C$ be a strongly $G$-graded tensor category. Then there is a bijective
correspondence between module categories over $\C_e$ and $G$-graded  $\C$-module categories.
\end{corol} \begin{proof} It  is a particular case of  Theorem \ref{estruc simples grad}, with $X=G$.
\end{proof} \begin{prop}\label{mapa picard} For every  $\sigma, \tau \in G$, the canonical functor
$$f_{\sigma,\tau}: \C_\sigma\boxtimes_{\C_e}\C_\tau\to \C_{\sigma\tau}, \ \ f_{\sigma,\tau}([X,Y])=
X\otimes Y,$$ is an equivalence of $\C_e$-bimodule categories. \end{prop}

\begin{proof} Let us consider the graded  $\C$-module category $\C(\tau)$, where $\C=\C(\tau)$ as
$\C$-module categories, but with grading $(\C(\tau))_g= \C_{\tau g}$, for $\tau\in G$.

Since $\C(\tau)_e=\C_\tau$, by Theorem \ref{estruc simples grad}, the canonical functor
$\mu(\C(\tau)_e): \C\boxtimes_{\C_\e}\C_\tau \to \C(\tau)$, \begin{equation*}[X,Y] \mapsto X\otimes Y
\end{equation*} is an equivalence of $G$-graded $\C$-module categories.  So the restriction
$\mu(\C(\tau))_\sigma:\C_\sigma \boxtimes_{\C_e}\C_\tau \to \C(\tau)_\sigma= \C_{\tau\sigma}$ is a
$\C_e$-module category equivalence. But by definition $\mu(\C(\tau))_\sigma= f_{\sigma,\tau}$. It is
clear that $f_{\sigma,\tau}$ is a $\C_e$-bimodule category functor, so the proof is finished.
\end{proof}

\section{Clifford Theory}

\emph{In this section we shall suppose that $\C$ is a strongly graded tensor category over a group
$G$.}

We shall denote by $\Omega_{\C_e}$ the set of equivalences classes of simple
$\C_e$-module categories. Given a $\C_e$-module category  $\M$, we shall denote by $\Omega_{\C_e}(\M)$
the set of equivalences classes of simple $\C_e$-submodule categories of $\M.$

\begin{lem}\label{accion indes}Let  $\M$ be a  $\C_e$-module category. Then for all $\sigma \in G$,
the category $\C_\sigma\boxtimes_{\C_e} \M$ is a simple $\C_e$-module category if and only if  $\M$
is. \end{lem} \begin{proof} If $\N$ is a proper $\C_e$-submodule category of $\M$, then the category
$\C_\sigma\boxtimes_{\C_e} \N$ is a $\C_e$-submodule category of   $\C_\sigma\boxtimes_{\C_e} \M$, so
the $\C_e$-module category $\C_\sigma\boxtimes_{\C_e} \M$ is not simple.

By Proposition \ref{mapa picard}, we have that $\M\cong
\C_{g^{-1}}\boxtimes_{\C_e}(\C_g\boxtimes_{\C_e}\M) $, so if $\C_g\boxtimes_{\C_e}\M$ is not simple,
then  $\M$ is not simple neither. \end{proof}

By Lemma \ref{accion indes} and  Proposition \ref{mapa picard}, the group   $G$ acts on
$\Omega_{\C_e}$ by $$G\times \Omega_{\C_e}\to \Omega_{\C_e}, \  \  (g,[X])\mapsto
[\C_g\boxtimes_{\C_e} X].$$

Let $\M$ be a $\C$-module category, and let $\N\subseteq \M$ be a Serre subcategory. We shall
denote by $\C_\sigma\overline{\otimes}\N$ the Serre  subcategory given by
$Ob(\C_\sigma\overline{\otimes}\N)=\{\text{subquotients of }$ $V\otimes N: V\in \C_\sigma, N\in \N\}$.
(Recall that a subquotient object is a subobject of a quotient object.)

\begin{prop}\label{lema g-conjunto} Let $\M$ be a  $\C$-module category and let $\N$ be a
$\C_e$-submodule category of $\M$. Then $\C_\sigma\boxtimes_{\C_e}\N \cong \C_\sigma\overline{\otimes}
\N$, as $\C_e$-module categories, for all $\sigma \in G$. \end{prop} \begin{proof} Define a $G$-graded
$\C$-module category by $gr\text{-}\N=\bigoplus_{\sigma\in G} \C_\sigma \overline{\otimes}\N$, with
action\begin{align*}
    \otimes: \C_\sigma \times \C_g \overline{\otimes}\N \to \C_{\sigma g}\overline{\otimes} \N\\
    V_\sigma\times T \mapsto V_\sigma\otimes T.
\end{align*} Since $\C_e\overline{\otimes} \N = \N$ as $\C_e$-module category, by Theorem \ref{estruc
simples grad} the canonical functor $\mu(\N):\C\boxtimes_{\C_e} \N \to gr-\N$ is a category
equivalence of $G$-graded  $\C$-module categories and the restriction  $\mu_\sigma:
\C_\sigma\boxtimes_{\C_e}\N\to \C_\sigma\overline{\otimes} \N$ is a $\C_e$-module category
equivalence. \end{proof}

\begin{corol}\label{G-conjuto} Let $\M$ be a $\C$-module category. The action of $G$ on
$\Omega_{\C_e}$ induces an  action of $G$ on $\Omega_{\C_e}(\M)$. \end{corol}

\begin{proof}

Let $\N$ be a simple  $\C_e$-submodule category of $\M$. By Proposition \ref{lema g-conjunto} the
functor \begin{align*}
    \mu_\sigma: \C_\sigma\boxtimes_{\C_e} \N&\to \C_\sigma\overline{\otimes} \N\\
                              [V, N] &\mapsto V\otimes N,
\end{align*}is a $\C_e$-module category equivalence, so  $\C_\sigma\boxtimes_{\C_e} \N$ is equivalent
to a  $\C_e$-submodule category of $\M$. \end{proof}

Let $\M$ be an abelian category and let  $\N, \N'$  be Serre subcategories of $\M$, we shall
denote $\N +\N'$ the Serre subcategory of $\M$ where $Ob(\N +\N')=\{ \text{subquotients  of } N
\oplus N' :N\in \N, N'\in \N' \}$. It will be called the sum category of $\N$ and $\N'$.

%Let $\M$ be a $\C$-module category and let  $\N$ be a  $\C$-submodule category of $\M$, then  \emph{the isotypic component category} of  $\N$ in $\M$ will be the sum category of all abelian module subcategories of $\M$ equivalent to $\N$ as $\C$-module categories. Now we are ready to give a proof of our main result.

\begin{proof}[Proof of the Theorem \ref{clifford}]

(1) \  Let $\N$ be a simple abelian $\C_e$-submodule category of $\M$, the canonical functor
\begin{align*}
    \mu: \C\boxtimes_{\C_e}\N&\to \M\\
                              [V,N] &\mapsto V\otimes N,
\end{align*} is a $\C$-module functor and $\mu=\oplus_{\sigma\in G} \mu_\sigma$, where
$\mu_{\sigma}=\mu|_{\C_\sigma}.$ By  Proposition \ref{lema g-conjunto} each $\mu_\sigma$ is  a
$\C_e$-module category equivalence with  $\C_\sigma\overline{\otimes}\N$.

Since $\M$ is simple, every object  $M\in \M$ is isomorphic to some subquotient of $\mu(X)$ for some object $X\in
\C\boxtimes_{\C_e}\N$. Then   $\M=\sum_{\sigma\in G}\C_\sigma\overline{\otimes} \N$ and each
$\C_\sigma\overline{\otimes} \N$ is an abelian simple  $\C_e$-submodule category.

Let $S, S'$ be simple abelian $\C_e$-submodule categories of $\M$. Then there exist  $\sigma, \tau \in
G$ such that $\C_\sigma\boxtimes_{\C_e}\N\cong S$, $\C_\tau\boxtimes_{\C_e}\N\cong S'$, and by
Proposition \ref{mapa picard}, $S'\cong \C_{\tau\sigma^{-1}}\boxtimes_{\C_e}S$. So the action is
transitive. \bigbreak

(2) Let  $H=\st([\N])$ be the stabilizer subgroup of  $[\N]\in \Omega_{\C_e}(\M)$ and let
$$\M_\N=\sum_{h\in H} \C_h\overline{\otimes} \N.$$

Since  $H$ acts transitively on $\Omega_{\C_e}(\M_\N)$, the $\C_H$-module category $\M_\N$ is simple.
Let $\Sigma= \{e,\sigma_1,\ldots\}$ be a set of representatives of the cosets of  $G$ modulo $H$. The
map $\phi: G/H\to \Omega_{\C_H}(M)$, $\phi(\sigma H)= [\C_\sigma\overline{\otimes}\M_\N]$ is an
isomorphism of $G$-sets. Then $\M$ has a structure of $G/H$-graded $\C$-module category, where $\M=$
$\oplus_{\sigma \in \Sigma}\C_\sigma\overline{\otimes} M_\N$. By Proposition \ref{estruc simples
grad}, $\M\cong \C\boxtimes_{\C_H}\M_\N$ as $\C$-module categories. \end{proof}

\begin{obs}
Nikshych and Gelaki noted the existence of a  grading by a transitive $G$-set for every indecomposable module category over a $G$-graded fusion category \cite[Proposition 5.1]{GeNi}. Using the Theorem \ref{estruc simples grad} and \cite[Proposition 5.1]{GeNi},  we can do an alternative proof of the main theorem in the case of strongly graded fusion categories.
\end{obs}

\section{Simple module categories over crossed product tensor categories and
$G$-equivariant of tensor categories} \subsection{$G$-equivariantization of tensor
categories} Let $G$ be a group acting on a  category (not necessarily by tensor
equivalences) $\C$,  $*: \underline{G}\to \underline{\text{Aut}(\C)},$ so we have the
following data \begin{itemize}
  \item functors $\sigma_*: \C\to \C$, for each $\sigma\in G$,
  \item isomorphism $\phi(\sigma,\tau): (\sigma\tau)_*\to \sigma_*\circ \tau_*$, for all $\sigma, \tau \in G$.
\end{itemize}   The category of $G$-invariant objects in $\C$, denoted by $\C^G$, is the category defined as  follows:  an object in $\C^G$ is a pair $(V, f)$, where $V$ is an object of $\mathcal M$ and $f$ is a family of isomorphisms $f_\sigma: \sigma_*(V) \to V$, $\sigma \in G$, such that, for all $\sigma, \tau \in G$, \begin{equation}\label{deltau} \phi(\sigma,\tau)f_{\sigma\tau}= f_\sigma \sigma_*(f_\tau).\end{equation} A $G$-equivariant morphism $\phi: (V, f) \to (W, g)$ between $G$-equivariant objects $(V, f)$ and $(W, g)$, is a morphism $u: V \to W$ in $\mathcal C$ such that $g_\sigma\circ \sigma_*(u) = u\circ f_\sigma$, for all $\sigma \in G$.

If the category $\C$ is a tensor category, and the action is by tensor
autoequivalences $*: \underline{G}\to \underline{\text{Aut}_\otimes(\C)},$ then we
have a natural isomorphism

\begin{itemize}
 \item  $\psi(\sigma)_{V,W}:\sigma_*(V)\otimes \sigma_*(W)\to \sigma_*(V\otimes
     W)$, for all $\sigma\in G$, $V, W\in \C$.
\end{itemize} thus $\C^G$ has a tensor product defined by \begin{align*}
    (V, f)\otimes (W, g):= (V\otimes W, h),
\end{align*}where $$h_\sigma= u_\sigma v_\sigma\psi(\sigma)_{V,W}^{-1},$$and unit
object $(1, \text{id}_1)$.

\begin{ejem}\textbf{The comodule category of a  cocentral cleft exact sequence of
Hopf algebras.} Let $G$ be a group and let \begin{equation}\label{extension Hopf}k\to
A\to H\overset{\pi}\rightarrow kG\to k \end{equation} be a cocentral cleft exact
sequence of Hopf algebras, \textit{i.e.}, the projection $\pi:H\to kG$ admits a
$kG$-colinear section $j:kG\to H$,  invertible with respect to convolution
product.

Since the sequence is cleft, the Hopf  algebra  $H$ has the structure of a bicrossed
product  $H\cong  A^\tau\#_\sigma kG$ with respect a certain compatible datum
$(\cdot, \rho, \sigma, \tau)$, where $\cdot : A\otimes kG\to A$ is a weak action,
$\sigma: kG\otimes kG\to A$ is invertible cocycle, $\rho: kG\to kG\otimes A$ is a
weak coaction, $\tau: kG\to A\otimes A$ is a dual cocycle, subjects to compatibility
conditions in \cite[Theorem 2.20]{AD}.

The projection in \eqref{extension Hopf}, is called cocentral if $\pi(h_1)\otimes h_2
=\pi(h_2)\otimes h_1$, this is equivalent to   the weak  coaction $\rho$ to be
trivial, see \cite[Lemma 3.3]{Na2}.

\begin{lem} Let $H\cong  A^\tau\#_\sigma kG$ be a   bicrossed product with trivial
coaction. Then the group $G$ acts over the category of right $A$-modules $\-_A\mathcal{M}$, and  $_H\mathcal{M}\cong (_A\mathcal{M})^G$ as tensor categories, were $_H\mathcal{M}$ is  the category of right $H$-modules . \end{lem}
\begin{proof} See \cite[Lemma 3.3]{Na2}. \end{proof} \begin{obs} Let $H$ be a
semisimple Hopf algebra  over $\mathbb C$. By \cite[Proof of Theorem 3.8]{GeNi}, the
fusion category $^H\M$ of finite dimensional comodules is $G$-graded (not necessary
strongly graded) if and only if there is a cocentral exact sequence of Hopf algebras
as in \eqref{extension Hopf}. In this case, the fusion category $^H\M$ is weakly
Morita equivalent to a $G$-crossed tensor category $_A\M\rtimes G$. That is,
$^H\M\cong \F_{_A\M\rtimes G}(\N,\N)$, for some indecomposable $_A\M\rtimes G$-module
category $\N$. \end{obs} \end{ejem}

\subsection{The obstruction to a $G$-action over a tensor category}

Let $\C$ be a tensor category, we shall denote  by $\text{Aut}_\otimes (\C)$  the
group of tensor auto-equivalences, it is the set of isomorphisms classes of
auto-equivalences of $\C$, with the multiplication induced by the composition:
$[F][F']= [F\circ F']$.

Every $G$-action over a tensor category induces a group homomorphism $\psi: G\to
\text{Aut}_\otimes (\C)$. We shall say that a homomorphism $\psi: G\to
\text{Aut}_\otimes (\C)$ is realizable if there is some $G$-action such the induced
group homomorphism coincides with $\psi$.

The goal of this subsection is show that for every homomorphism $\psi: G\to
\text{Aut}_\otimes (\C)$,  there is an associated element in a 3rd cohomology group
which is zero if and only if $\psi$ is realizable. Moreover, every realization is in
correspondence (non natural) with an element of a 2nd cohomology group.

\subsubsection{Categorical-groups}\label{cat-group} A categorical-group $\mathcal{G}$ is a monoidal category where every object, and every arrow is invertible, see \cite{Baez} for a complete reference.

A trivial example of a categorical-group is the discrete categorical-group $\underline{G}$, associated to a group $G$. The objects of $\underline{G}$ are the elements of $G$, the arrows are only the identities, and the tensor product is the multiplication of $G$.

Complete invariants of a categorical-group $\mathcal{G}$ with respect to monoidal equivalences
are
$$\pi_0(\mathcal{G}), \pi_1(\mathcal{G}), \alpha\in H^3(\pi_0(\mathcal{G}),\pi_1(\mathcal{G})),$$ where $\pi_0(\mathcal{G})$ is the group of isomorphism classes of objects, $\pi_1(\mathcal{G})$ is the abelian group of automorphisms of the unit object. The group $\pi_1(\mathcal{G})$ is a $\pi_0(\mathcal{G})$-module in the natural way, and $\alpha$ is a third cohomology class given by the associator.

Complete invariants of a monoidal functor $F: \mathcal{G}\to \mathcal{G}'$ between categorical-groups, with respect to monoidal isomorphisms are

$$\pi_0(F):\pi_0(\mathcal{G})\to \pi_0(\mathcal{G}'),  \pi_1(F): \pi_1(\mathcal{G})\to \pi_1(\mathcal{G}'), \theta(F):\pi_0(\mathcal{G})\times \pi_0(\mathcal{G})\to \pi_1(\mathcal{G}')$$ where $\pi_0(F)$ is a morphism of groups, $\pi_1(F)$ is a morphism of $\pi_0(\mathcal{G})$-modules and $\theta(F)$ is a class in $C^2(\pi_0(\mathcal{G}),\pi_1(\mathcal{G}'))/B^2(\pi_0(\mathcal{G}),\pi_1(\mathcal{G}'))$, such that $$\delta(\theta(F))=\pi_1(\mathcal{G}')_*(\phi(\mathcal{G}))-\pi_0(\mathcal{G}')^*(\phi(\mathcal{G}')),$$ where
\begin{align*}
    \pi_0(F)^*:C^*(\pi_0(\mathcal{G}'),\pi_1(\mathcal{G}'))&\to C^*(\pi_0(\mathcal{G}),\pi_1(\mathcal{G}')),\\
    \pi_1(F)_*:C^*(\pi_0(\mathcal{G}),\pi_(\mathcal{G})) &\to C^*(\pi_0(\mathcal{G}),\pi_1(\mathcal{G}')),
\end{align*}are the maps of cochain complexes induced by the group morphisms $\pi_0(F)$ and $\pi_1(F)$.
The next result follows from the last discussion or see \cite{Baez}.
\begin{prop}\label{prop cat-groups}
Let $\mathcal{G}$ be a categorical group and let $f:G\to \pi_0(\mathcal{G})$ be a morphism of groups. Then there is a monoidal functor $F:\underline{G}\to \mathcal{G}$, such that $f=\pi_0(F)$ if and only if the cohomology class of $f_*(\phi)$ is zero.

If  $f_*(\phi)$ is zero, the classes  of equivalence of monoidal functors $F:\underline{G}\to \mathcal{G}$ are in one to one correspondence with $H^2(G,\pi_1(\mathcal(\G)))$.
\end{prop}
\begin{proof}
The monoidal category $\underline{G}$ has invariants $\pi_0(\underline{G})= G$ and $\pi_1(\underline{G})=0$. Then, the proof follows from the discussions of this subsection, or see \cite{Baez}.
\end{proof}
\subsubsection{The obstruction to a $G$-action over a tensor category and cyclic
actions}

Let  $\underline{\text{Aut}_\otimes(\C)}$ be the monoidal category of tensor
auto-equivalences of a tensor category $\C$, where arrows are tensor natural
isomorphisms and tensor product given by composition of functors. Then
$\underline{\text{Aut}_\otimes(\C)}$ is a categorical-group.

The invariants associated to $\underline{\text{Aut}_\otimes(\C)}$ (see Subsection \ref{cat-group}) are  Then $\pi_0(\underline{\text{Aut}_\otimes(\C)})= \text{Aut}_\otimes (\C)$, and $\pi_1(\underline{\text{Aut}_\otimes(\C)}) = \text{Aut}_\otimes(\text{id}_\C)$, the group of monoidal natural isomorphisms of the identity functor.

\begin{teor}\label{obstruccion} Let $\C$ be a tensor category and let $G$ be a group.
Consider the data $(\text{Aut}_\otimes (\C),\text{Aut}_\otimes(\text{id}_\C),[a])$
associated to the categorical-group $\underline{\text{Aut}_\otimes(\C)}$. Then \begin{itemize}
\item a group homomorphism $f: G\to \text{Aut}_\otimes (\C)$ is realized as a
$G$-action over $\C$ if and only if $0=[f_*(a)]\in
H^3(G,\text{Aut}_\otimes(\text{id}_\C))$.

\item If the group homomorphism $f: G\to \text{Aut}_\otimes (\C)$ is realizable,
    then the set of realizations of $f$ is in 1-1 correspondence with
    $Z^2(G,\text{Aut}_\otimes(\text{id}_\C))$,  and the set of equivalences classes
    of realizations of $f$ is in 1-1 correspondence with
    $H^2(G,\text{Aut}_\otimes(\text{id}_\C))$. \end{itemize} \end{teor} \begin{proof}
    The Theorem is a particular case of the Proposition \ref{prop cat-groups}. \end{proof}

Recall  that if $A$ is a module for the cyclic group $C_m$ of order $m$, then:
\begin{align}\label{cohomologia ciclico} H^n(C_m; A)= \begin{cases}\{a\in A:
Na=0\}/(\sigma-1)A, \qquad &\text{if } n=1, 3, 5, \ldots \\ A^{C_m}/NA, \quad
&\text{ if } n = 2, 4, 6, \ldots, \end{cases} \end{align} where $N= 1 + \sigma
+\sigma^2 +\cdots + \sigma ^{m-1}$, see \cite[Theorem 6.2.2]{Weibel}. Given an
element $a \in A^{C_m}$ the associated 2-cocycle  can be constructed as follows.
\begin{align}\label{2-cociclo ciclico} \gamma_a(\sigma^i, \sigma^j)= \begin{cases}1,
\quad &\text{if } i+j < m, \\ a^{i+j-m}, \quad  &\text{ if } i+j\geq m. \end{cases}
\end{align}

Let $F: \C\to \C$ be a monoidal equivalence, such that there is a monoidal natural
isomorphism $\alpha: F^m\to  \text{id}_\C$. By Theorems \ref{obstruccion} and
\eqref{cohomologia ciclico}, the induced homomorphism $\psi: C_m\to
\text{Aut}_\otimes (\C)$ is realizable if and only if $\text{id}_F\otimes
\alpha\otimes \text{id}_{F^{-1}}= \alpha$. In this case, two natural isomorphisms
$\alpha_1, \alpha_2: F^m\to \text{id}_\C$ realize equivalent $C_m$-actions if and
only if there is a monoidal natural isomorphism $\theta: F_1\to F_2$ such that
$\theta^mF_1=F_2$.

\begin{corol}\label{cicly actions} Let $\C$ be a tensor category and let $C_m$ be
cyclic group of order $m$. Then the set of $C_m$-actions over  $\C$ are in 1-1
correspondence with pairs $(F, \alpha)$, where $F:\C \to \C$ is a monoidal
equivalence, $\alpha:F^m\to \text{id}_\C$ is a monoidal natural isomorphism such
$\text{id}_F\otimes \alpha =\alpha\otimes \text{id}_F$.

Two pairs $(F_1,\alpha_1)$ and $(F_2,\alpha_2)$ induce equivalent $C_m$-actions if
and only if there is a monoidal natural isomorphism $\theta: F_1\to F_2$ such that
$\theta^mF_1=F_2$.  \end{corol}

The description of the 2-cocycle associated to a $C_m$-invariant element \eqref{2-cociclo ciclico}, is as follows: the $C_m$-action $\psi: C^m\to \underline{\text{Aut}_\otimes(\C)}$
associated to a pair $(F,\alpha)$ is $\psi(1)= \text{id}_\C$, $\psi(\sigma^i)= F^i$,
$i=1,\ldots m-1$, and the monoidal natural isomorphisms
$\phi_\alpha(\sigma^i,\sigma^j): F^i\circ F^j\to F^{i+j}$

\begin{align}\phi_\alpha(\sigma^i, \sigma^j)= \begin{cases}\text{id}_\C, \quad
&\text{if } i+j < m, \\ \text{id}_F\otimes\alpha^{i+j-m}= \alpha^{i+j-m}\otimes
\text{id}_F, \quad  &\text{ if } i+j\geq m. \end{cases} \end{align}

\subsubsection{The bigalois group of a Hopf algebra}

Let $H$ be a Hopf algebra. A right $H$-Galois object is a non-zero right $H$-comodule
algebra $A$ such that the linear map defined by $\text{can}: A\otimes A\to A\otimes
H, a\otimes b\mapsto ab_{(0)}\otimes b_{(1)}$ is bijective.

A fiber functor $F:\ ^H\M\to Vec_k$  is an exact and faithful monoidal functor that
commutes with colimits. Ulbrich defined in \cite{Ulbrich}  a fiber functor $F_A$
associated with each $H$-Galois object $A$, in the form  $F_A(V)= A\square_H V$,
where $A\square_H V$ is the cotensor product over $H$ of the right $H$-comodule $A$
and the left $H$-comodule $V$. He showed in \textit{loc. cit.} that this defines a
category equivalence between $H$-Galois objects and fiber functors over $^H\M$.

 Similarly, a left $H$-Galois object is a non-zero left $H$-comodule algebra
$A$ such that the linear map $\text{can}: A\otimes A\to H\otimes A, a\otimes b
\mapsto a_{(-1)}\otimes a_{(0)}b$ is bijective.

Let $H$ and $Q$ be Hopf algebras.  An $H$-$Q$-bigalois object is an algebra $A$ which
is an $H$-$Q$-bicomodule algebra and both a left $H$-Galois object and a right
$Q$-Galois object.

Let $A$ be an $H$-Galois object. Schauenburg shows in \cite[Theorem 3.5]{scha-biga
communi alg} that there is a Hopf algebra $L(A,H)$ such that $A$ is a
$L(A,H)$-$H$-bigalois object.

The Hopf algebra $L(A,H)$ is the Tannakian-Krein reconstruction from the fiber
functor associated to $A$. By   \cite[Corollary 5.7]{scha-biga communi alg}, the
following categories are equivalent: \begin{itemize} \item The monoidal category
$\underline{\text{BiGal}(H)}$, where objects are $H$-bigalois object, morphism are
morphism of $A$--bicomodules algebras, and tensor product $A\square_H B$, the
cotensor product over $H$. \item The monoidal category
$\underline{\text{Aut}_\otimes(\ ^H\M)}$. \end{itemize}

Schauenburg defined the group $\text{BiGal}(H)$ as  the set of isomorphism classes of
$H$-bigalois objects  with multiplication induced by the cotensor product. This group
coincides with $\text{Aut}_\otimes (^H\M)$.

Is easy to see that for the Hopf algebra $kG$ of a group $G$, $\text{BiGal}(kG)=
Aut(G)\rtimes H^2(G,k^*)$. However, it is difficult to find an explicit description
in general. The group $\text{BiGal}(H)$ has been calculated for some Hopf algebras,
for example: Taft algebras \cite{scha-taft}, monoidal non-semisimple Hopf algebras
\cite{Bichon}, the algebra of function over a finite group coprime to $6$
\cite{davydov}.

\subsubsection{The abelian group $\text{Aut}_\otimes(\text{id}_{\C})$ for Hopf
algebras}

\begin{prop}\label{group-like centrales} Let $H$ be a Hopf algebra. Then
$\text{Aut}_\otimes(\text{id}_{_H\M})\cong G(H)\cap Z(H)$ the group of central
group-likes of $H$. \end{prop}

\begin{proof} The maps $H\otimes_k(M\otimes_k N)\to (H\otimes_H
N)\otimes_k(H\otimes_H N), h\otimes m\otimes n \mapsto (h_{(1)}\otimes m)\otimes
(h_{(2)}\otimes n)$, and $H\otimes_k k\to k, h\otimes 1\mapsto \epsilon(h)$, induce
natural $H$-module morphisms \begin{align*}
    F_{M,N}:H\otimes_H(M\otimes N)&\to (H\otimes_H M)\otimes (H\otimes_H N)\\
    F^0: H\otimes_H k &\to k.
\end{align*} The identity monoidal functor is naturally isomorphic  to $(_\cdot
H_\cdot \otimes_H (-), F, F^0)$, and it is well-know that every $H$-bimodule
endomorphism is of the form $\psi_c:H\to H, h\mapsto ch$, for some $c\in Z(H)$. The
natural transformation associated to $\psi_c$ is monoidal if and only if $\psi_c$ is
a bimodule coalgebra map, \emph{i.e.}, if $c$ is a group-like. \end{proof}

For the group algebra $kG$, we have $\text{Aut}_\otimes(\text{id}_{_{kG}\M})\cong
Z(G)$ the center of $G$, and for a Hopf algebra $\mathbb C^G$, where $G$ is a finite
group, we have $\text{Aut}_\otimes(\text{id}_{_{k^G}\M})\cong G/[G,G]$.

 Let $\C$ be a complex fusion category, $i.e.$, a semisimple tensor category
with finitely many isomorphisms classes of simple objects. In \cite{GeNi} it is shown
that every fusion category is naturally graded by a group $U(\C)$  called the
universal grading group of $\C$. The group $U(\C)$ only depends of the Grothendieck
ring of $\C$.

In \cite[Proposition 3.9]{GeNi} it is shown  that if $\C$ is a fusion category and
$G= U(\C)$ is the universal grading group of $\C$, then
$\text{Aut}_\otimes(\text{id}_{\C})\cong \widehat{G_{ab}}$ the group of characters of
the maximal abelian quotient of $G$. \begin{corol} Let $H$ be a semisimple
almost-cocommutative Hopf algebra. Then $U(_H\M)\cong Z(H)\cap G(H)$. \end{corol}

\begin{proof} Since $H$ is almost-commutative the Grothendiek ring is commutative,
hence the universal grading group is abelian. By Proposition \ref{group-like
centrales} and \cite[Proposition 3.9]{GeNi} $U(_H\M)\cong Z(H)\cap G(H)$.
\end{proof}

\subsection{$G$-invariant actions on module categories}

Let $\C$ be a tensor category and let $(\sigma, \psi): \C\to \C$ be a monoidal functor. If $(\M,\otimes, \alpha)$ is a right $\C$-module category, the twisted  $\C$-module category $(\M^\sigma, \otimes^\sigma, \alpha^\sigma)$ is defined by: $\M=\M^\sigma$ as category, with $M\otimes^\sigma V = M\otimes \sigma (V)$, and $\alpha^\sigma_{M,V,W}=\text{id}_M\otimes \psi_{V,W}\circ\alpha_{M,V,W}$.  %If we define the $\C$-bimodule $\widehat{\C}^\sigma$, where $\widehat{\C}^\sigma =\C$ as left $\C$-module category, $\widehat{\C}^\sigma = \C^\sigma$ as right $\C$-module category and associativity constraint $$\alpha^\sigma_{V,W,Z}= \alpha_{V,W,\sigma_*(Z)},$$ for all $V, Z\in \C, W\in \C^\sigma$. Then the canonical $\C$-functor $\widehat{\C}^\sigma\boxtimes_{\C} \M\to \M^\sigma, [V,M]\mapsto V\otimes M$ is an equivalence of $\C$-module categories.

\begin{defin} Let $\C$ be a tensor category,  $\M$  a left $\C$-module category and
$\sigma:\C\to \C$ a monoidal functor. We shall say that the functor $(T,\eta):M\to
M^{\sigma}$ is a $\sigma$-equivariant functor of $\M$ if is a $\C$-module functor.

Given an action of a group  $G$ over $\C$, the module category $\M$ is called
$G$-invariant if there is  a $\sigma$-equivariant functor for each $\sigma\in G$.
\end{defin}

Let $\sigma, \tau:\C\to \C$ be monoidal functors. Let also $(T,\eta):M\to M^{\sigma}$
a $\sigma$-equivariant functor and $(T',\eta'):M\to M^{\tau}$ a $\tau$-invariant
functor. We define their composition by  $$(T'T,T'(\eta)(\eta'(T\times T))):\M\to
\M.$$ This gives a $\sigma\circ \tau$-equivariant functor of $\M$.

 Given a $G$-action over a monoidal category $\C$  and a $G$-invariant
module category $\M$, we denote by $\underline{\text{Aut}^G_{\C}(\M)}$ the following
monoidal category: objects are $\sigma_*$-equivariant functors, for all $\sigma \in
G$, morphisms are natural isomorphisms of module functors, the tensor product is
composition of $\C$-module functors and the unit object is the identity functor of
$\M$.

\begin{defin} Let $(\sigma_*,\phi(\sigma,\tau),\psi(\sigma)):\underline{G} \to
\underline{\text{Aut}_\otimes (\C)}$ be  an action of $G$ over a tensor category
$\C$, and let $\M$ be a $G$-invariant $\C$-module category. A $G$-invariant functor
over $\M$ is a monoidal functor $(\sigma^*, \phi,\psi):\underline{G} \to
\underline{\text{Aut}^G_{\C}(\M)}$, such that $\sigma^*$ is a $\sigma_*$-invariant
functor, for all $\sigma\in G$. \end{defin}

\begin{obs}(1) A $\C$-module category $\M$ with a $G$-invariant functor is called a
\emph{$G$-equivariant $\C$-module category} in \cite[definition 5.2]{ENO2}.

(2) Let $\C$ be a $G$-invariant monoidal category. The monoidal category
$\underline{\text{Aut}^G_{\C}(\M)}$ is a graded categorical-group and the group
$\text{Aut}^G_{\C}(\M)$ has a natural group epimorphism $\pi:
\text{Aut}^G_{\C}(\M)\to G$. So, if a group homomorphism $\psi:G\to
\text{Aut}^G_{\C}(\M)$ is realizable, then $\pi \psi= \text{id}_G$. Such group
homomorphisms will be called \emph{split}.

(3) Let $\psi:G\to \text{Aut}^G_{\C}(\M)$ be a split group homomorphism. If  $a\in
H^3$ $(\text{Aut}^G_{\C}(\M), H)$ is the 3-cocycle associated to the categorical-group
$\underline{\text{Aut}^G_{\C}(\M)}$, then  like in Theorem \ref{obstruccion}, $\psi$
is realizable if and only if the 3-cocycle $\psi_*(a)$ is a 3-coboundary, and the set
of realizations of $\psi$ is in correspondence with the elements of a 2nd cohomology
group.

\end{obs} The following result appears in \cite[Sec. 2]{tambara}.
\begin{prop}\label{funtor invariante} Let $\C\rtimes G$ be a crossed product tensor
category. Then there is a bijective correspondence between structures of $\C\rtimes
G$-module category  and $G$-invariant functors over a $\C$-module category $\M$.
\end{prop}

\begin{proof} Let $\M$ be a $\C\rtimes G$-module category. Each object $[1, \sigma]$,
$\sigma\in G$, defines an equivalence $\sigma_* :\M\to \M, M\mapsto [1,\sigma]\otimes
M$. With $\phi(\sigma,\tau)_M= \alpha_{(1,\sigma),(1,\tau),M}$ the constraint of
associativity, this defines a monoidal functor $\underline{G}\to \text{Aut}(\M)$.

The category  $\M$ is a $\C$-module category with $V\otimes M = [V,e]\otimes V$ and
since $[1,\sigma]\otimes [V,e]= [\sigma_*(V),e]\otimes [1,\sigma]$ we have a natural
isomorphism $\psi(\sigma)_{V,M}: \sigma_*(V)\otimes \sigma (M)\to \sigma (V\otimes
M)$, by $\psi(\sigma)_{V,M}= \alpha_{(1,\sigma),(V,e),M}^{-1}\circ
\alpha_{(\sigma_*(V),e),(1,\sigma),M}$. This defines a $G$-invariant functor.

Conversely, if $\underline{G} \to \text{Aut}^G_{\C}(\M)$ is a $G$-invariant functor,
we have natural isomorphisms $\phi(\sigma,\tau)_M:\sigma_*\tau_*(M)\to
\sigma\tau_*(M)$, $\psi(\sigma)_{V,M}:\sigma(V)\otimes \sigma (M)\to \sigma(V\otimes
M)$. Then, we may define the action on $\M$ by $$(V,\sigma)\otimes M:= V\otimes
\sigma_*(M),$$ and constraint of associativity $$\alpha_{(V,\sigma), (W,\tau),M} =
\text{id}_{V\otimes \sigma_*(W)} \otimes \phi(\sigma,\tau)_M\circ
\alpha_{V,\sigma_*(W),\sigma_*(\tau_*(M))} \circ \text{id}_V\otimes
\psi(\sigma)_{W,M}^{-1}.$$ \end{proof}

Suppose that the group $G$ is finite and the tensor category $\C$ is a fusion
category over an algebraically closed field of characteristic zero. Then the module
categories over $\C\rtimes G$ and $\C^G$ are in bijective correspondence by
\cite[Proposition 3.2]{nik}. If $\M$ is $\C\rtimes G$-module category then, by
Proposition \ref{funtor invariante}, there is a $G$-action on $\M$, and the category
$\M^G$ is a $\C^G$-module category with \begin{align*}
    (V, f)\otimes (M, g):= (V\otimes M, h),
\end{align*}where $$h_\sigma= g_\sigma h_\sigma\psi(\sigma)_{V,M}^{-1}.$$ For a
$k$-linear monoidal category and $G$ finite where char$(k)\nshortmid |G|$, Theorem
\cite[Theorem 4.1]{tambara} says that every $\C^G$-module category is of the form
$\M^G$ for a $\C\rtimes G$-module category. The following result appears in
\cite{ENO2} for fusion categories and finite groups.

\begin{teor}\label{clasi module  crossed} Simple module categories over $\C\rtimes G$
are in bijective corres\-pondence with the following data:

\begin{itemize} \item  a subgroup $H\subseteq G$,
  \item a simple $H$-invariant $\C$ module category $\M$,
  \item a monoidal functor $\underline{H} \to \text{Aut}^H_{\C}(\M)$. %$(\sigma_*:\phi(\sigma, \tau),\nu)\underline{H}\to \underline{\text{Aut}(\M)}$
\end{itemize} If the group $G$ is finite then the module categories over $\C^G$
are in bijection with the same data.

\end{teor} \begin{proof} By Theorem \ref{clifford}, if $\N$ is an simple $\C\rtimes
G$-module category, then  it is isomorphic to $\C\boxtimes_{\C\rtimes H}\M$ for some
subgroup $H\subseteq G$ and a simple  $\C\rtimes H$-module category $\M$, such that
$\M$ is $H$-invariant. In particular it follows that the restriction of $\M$ to $\C$
is simple. Now the correspondence follows from Proposition \ref{funtor invariante}.

If the group $G$ is finite, then the  correspondence follows from \cite[Theorem
4.1]{tambara} or \cite[Proposition 3.2]{nik}. \end{proof}

Suppose that $G$ is a finite group and   $H\cong  A^\tau\#_\sigma kG$ is a bicrossed
product with trivial coaction. Then the module categories over $_H\mathcal{M}$ are of
the form $\N^G$, for some $G$-equivariant $_A\M$-module category $\N$. Moreover, the
module category is simple if and only if $\N$ is simple.

\begin{ejem} Let $N \geq 2$ be an integer and let  $q \in \mathbb C$  be a primitive
$N$-th root of unity.  The {\it Taft algebra} $T(q)$ is the $\mathbb C$-algebra
presented by generators $g$ and $x$ with relations $g^N = 1$, $x^N = 0$ and $gx = q
xg$. The algebra $T(q)$ carries a Hopf algebra structure, determined by $$ \Delta g =
g \otimes g, \qquad \Delta x = x \otimes 1 + g \otimes x.$$ \noindent Then $
\varepsilon (g) = 1,\ \varepsilon(x) = 0,\ \mathcal S(g) = g^{-1}$, and $\mathcal
S(x) = -g^{-1}x$. It is known that \begin{enumerate} \item $T(q)$ is a pointed
non-semisimple Hopf algebra, \item the group of group-like elements of $T(q)$ is
$G(T(q)) = \langle g  \rangle \simeq {\mathbb Z} / (N)$, \item $T(q) \simeq T(q)^*$,
\item $T(q) \simeq T(q')$ if and only if $q = q'$. \end{enumerate}

\begin{prop}\label{taft} Let $G$ be a group, then the set of $G$-actions on  the
tensor category $\ ^{T(q)}\M$ of $T(q)$-comodules  is in 1-1 correspondence with the
set of group homomorphism from $G$ to $\mathbb C^*\ltimes \mathbb C$, where $\mathbb
C^*$ acts on $\mathbb C$ by $\mathbb C^* \times \mathbb C\to \mathbb C, (s,t)\mapsto
st$. \end{prop}

\begin{proof} By Proposition \ref{group-like centrales},  the abelian group
$\text{Aut}_\otimes(\text{id}_\C)$ is trivial, and by \cite[Theorem 5]{scha-taft},
$\text{Aut}_\otimes (^{T(q)}\M)= \text{BiGal}(T(q))\cong \mathbb C^*\rtimes \mathbb
C$. Then by Theorem \ref{obstruccion}, the set of isomorphism classes of $G$-actions
is given by the set of group homomorphism from $G$ to $\mathbb C^*\ltimes \mathbb C$.
\end{proof}

If $G=\mathbb Z/ (N)$ then, by Proposition \ref{taft} the possible $G$-actions are
parameterized by pairs $(r,\mu)$, where $r$ is a non-trivial  $N$-th root of the unit
and  $\mu\in \mathbb C$.

We shall denote by $A_{(\alpha,\gamma)}$ the $T(q)$-bigalois object associated to the
pair $(r,\mu)\in  \mathbb C^*\ltimes \mathbb C\cong \text{BiGal}(T(q))$. See
\cite[Theorem 5]{scha-taft}.

The ${}^{T(q)}\M$-module categories of rank one are in correspondence with fiber
functors on ${}^{T(q)}\M$, and these are in turn in 1-1 correspondence with
$T(q)$-Galois objects. By Theorem 2 in \emph{loc. cit.}, every $T(q)$-Galois object
is isomorphic to $A_{(1,\beta)}$, $\beta\in\mathbb C$, and two $T(q)$-Galois objects
$A_{(1,\beta)},$ $A_{(1,\mu)}$ are isomorphic if and only $\beta=\mu$.

 By Theorem \ref{clasi module  crossed}, if there is a semisimple module
category of rank one  over $\C =\ ^{T(q)}\M\rtimes \mathbb Z/ (N)$, it must be a
$^{T(q)}\M$-module category $\mathbb Z/ (N)$-invariant.

Suppose that $A_{(1,\beta)}$ is $\mathbb Z/(N)$-invariant. Since
$A_{(r,\mu)}\square_{T(q)} A_{(1,\beta)}\cong A_{(r,\mu+\beta)}$, we have that
$\mu=0$. Then if the action is associated to a pair $(r,\mu)$ where $\mu \neq 0$, the
category $\C$ does not admit any fiber functor, \emph{i.e.}, it is not the category
of comodules of a Hopf algebra.

However, since every simple object is invertible, the Perron-Frobenius dimension of
the simple objects is one. So, by \cite[Proposition 2.7]{finite-categories}, the
tensor category $\ ^{T(q)}\M\rtimes \mathbb Z/ (N)$ is equivalent to the category of
representations of a quasi-Hopf algebra.

 Note that the tensor category $(^{T(q)}\M)^{G}$ has at least one fiber
functor, for every group and every group action. In fact, since the forgetful functor
$U:\ ^{T(q)}\M^{G}\to\ ^{T(q)}\M$ is monoidal, then the composition with the fiber
functor of $^{T(q)}\M$ gives a fiber functor on $(^{T(q)}\M)^{G}$. \end{ejem}

\medbreak \noindent \textbf{Acknowledgement.} 
The author thanks Pavel Etingof, Gast\'on Garc\'ia, Mart\'in Mombelli,   Sonia Natale, and Dmitri Nikshych for useful discussions and advice.

\end{document}